\documentclass[12pt]{article}
\usepackage{amsmath,amsfonts,amssymb,amscd}
\usepackage{verbatim}
\title{Quartic double solids with ordinary singularities}
\author{M.I. Grooten \and J.H.M. Steenbrink}
\date{30 June 2008}
\newcommand{\PP}{{\mathbb{P}}}

\newcommand{\C}{{\mathbb{C}}}
\newcommand{\Q}{{\mathbb{Q}}}

\newcommand{\Z}{{\mathbb{Z}}}

\newcommand{\calO}{{\cal O}}

\newcommand{\calD}{{\cal D}}

\newtheorem{theorem}{Theorem}
\newtheorem{corollary}[theorem]{Corollary}
\newtheorem{lemma}[theorem]{Lemma}

\usepackage{hyperref}

\begin{document}
\maketitle
\begin{abstract}
We study the mixed Hodge structure on the third homology group of a threefold which is the double cover of projective three-space ramified over a quartic surface with a double conic. We deal with the Torelli problem for such threefolds.
\end{abstract}
\section{Introduction} \emph{Double solids}, that are: double covers of the projective
3-space, have been studied by, among others, Welters \cite{W} and Clemens \cite{DS}.
The latter admitted ordinary double points in his double solids. Our aim 
is to study double solids whose singularities are
1-dimensional, in particular those whose singularities are of type
$A$, $D$ or $T$; these varieties are called \emph{ordinary double
solids}.

A property of ADT-threefolds, which is a more general class of
threefolds, that includes the ordinary double solids, is that the
Hodge structures on all of their homology groups are pure, except
that on $H_3$, which might have a non-empty part of weight two
(Theorem \ref{ADT}). Moreover, this mixed Hodge structure can be considered as an extension of pure Hodge
structures, and the extension class  is strongly related to the Abel-Jacobi map for
$H_3(\hat{X})$, where $\hat{X}$ is a resolution of the ADT threefold
$X$ (Theorem \ref{extensionclass}). 

Non-trivial ordinary double solids first occur in degree four, where
there are 11 classes. The varieties in two of these classes fail to have a
pure Hodge structure on $H_3$. The \emph{cyclide double solids},
defined as double covers of $\PP3(\C)$ ramified along a quartic
surface singular along a plane conic, are one class. In the final two
paragraphs of this article, their full mixed Hodge structures are
computed. We conclude that the functor $X \mapsto H_3(X,\C)$, which
maps a cyclide double solid to the polarized mixed Hodge structure on
its third homology group, is six-to-one and is characterized by
'forgetting' the choice of one point among a six-tuple of points on
$\PP1(\C)$.

We intend to treat the other interesting class, the \emph{spider double solids} in a future paper. 
\section{ADT threefolds}
We first recall (cf. \cite[pp.~616--618]{GH})
the notion of ordinary singularities for surfaces. Ordinary surface singularities are the singularities arising as singular points
of the images of generic projections of smooth complex projective surfaces into  $\PP^3$. There are three types of these singularities:
\begin{enumerate} 
\item A \emph{transverse double point} with local equation $uv=0$ in an appropriate local holomorphic coordinate system.
\item A \emph{triple point} with local normal form $uvw=0$.
\item A \emph{pinch point} with local normal form $u^2-vw^2$.
\end{enumerate}

In this paper we will deal with threefolds which are double covers of $\PP^3$. These are called \emph{double solids}. The singularities of a double solids are exactly the preimages of the singular points of its ramification surface. We call a singular point of a double solid \emph{ordinary} if the corresponding singularity of its ramification surface is ordinary. Hence we obtain the following three types of ordinary singularities of double solids:
\begin{description}
\item[Type $A$] with local normal form $t^2=uv$.
\item[Type $T$] with local normal form $t^2=uvw$.
\item[Type $D$] with local normal form $t^2=u^2-vw^2$.
\end{description}
An \emph{ADT-threefold} is by definition a three-dimensional complex projective variety with at most singularities of types A, D or T. 

\begin{lemma} The points of type $A$ and $T$ are quotient singularities. 
\end{lemma}
\paragraph{Proof} For type $A$ consider the involution $(x,y,w) \mapsto (-x,-y,w)$ of $\C^3$ with invariants generated by $w,t=xy,u=x^2,v=y^2$. For type $T$ consider the group of diagonal matrices $\mathrm{diag}(a,b,c)$ with $a^2=b^2=c^2=abc=1$ acting on $\C^3$ with coordinates $x,y,z$. The invariant ring is generated by $t=xyz,u=x^2,v=y^2,w=z^2$ and $t^2=uvw$. 
\begin{theorem}
Let $X$ be an ADT-threefold with reduced singular locus $\Sigma$. Then the blowing up of $\Sigma$ in $X$ is a resolution $\sigma:\hat{X} \to X$. Moreover $\sigma^{-1}(\Sigma)=E$ is a divisor with normal crossings on $\hat{X}$. The fibres of $\sigma$ over points of $\Sigma$ are: 
\begin{enumerate}
\item $\PP^1$ for a point of type  $A$;
\item two intersecting $\PP^1$'s for a point of type $D$;
\item three  $\PP^1$'s  intersecting in one point for a point of type $T$.
\end{enumerate}
\end{theorem}
The proof is elementary; one just calculates starting with the three normal forms. Let us treat the case of type $T$, and consider the normal form $t^2=uvw$. The singular locus $\Sigma$ is the union of the $u$-,$v$- and $w$-axes,
and its ideal is  $I(\Sigma)=(t,uv,uw,vw)$. Over the affine open subset $t\neq 0$, the blowing up map $\sigma:\hat{X} \to X$ is given by 
$$\sigma(x,y,z)=(xyz,yz,xz,xy)\ .$$
Over the affine subset where $uv\neq 0$ which by symmetry represents the two remaining cases as well, we have 
$$\sigma(x,y,z)=(xyz^2,x,y,xyz)\ .$$
\begin{theorem} \label{ADT}
Let $X$ be an ADT-threefold. Define its \emph{defect} as $\mathrm{def}(X):= b_4(X)-b_2(X)$ where $b_i(X)$ denotes the i-th Betti number (cf. \cite[Theorem 3.2]{NS} and \cite[Sec.~3]{DS}). Then 
\begin{enumerate}
\item the mixed Hodge structure on $H^k(X)$ is pure for all $k\neq 3$;
\item $W_1H^3(X)=0$ and $W_2H^3(X)$ is purely of type $(1,1)$;
\item $$\mathrm{def}(X) = \delta - \dim W_2H^3(X)$$ where $\delta$ is the number of singular points on $X$ of type $D$. 
\end{enumerate}
\end{theorem}
\paragraph{Proof} Let $U \subset X$ denote the complement of the set $\calD$ of singular points of type $D$. Note that $H^i_{\calD}(X,\Q)=0$ unless $i=3,4$ or $6$ and that $H^3_{\calD}(X,\Q)$ is a Hodge structure which is purely of type $(1,1)$ and of rank the cardinality of $\calD$. 

To see this note that the singularity of type $D$ covers and is covered by the ordinary double point in dimension three. To make this explicit, change coordinates in the normal form $t^2=u^2-vw^2$. Put $p=u+t$ and $q=u-t$ so $xy=uv^2$. Then $P=p^2,\ Q=q^2,\ V=v^2$ and $W=w$ satify $PQ=VW$. In this way the singularity of type $D$ is sandwiched between two copies of the ordinary double point mapped one to another by $(P,Q,V,W) \mapsto (P^2,Q^2,V^2,W^2)$. Because this endomorphism of the ordinary double point induces an isomorphism on its rational local cohomology groups, the singularity of type $D$ must have its local cohomology groups with rational coefficients isomorphic to those of the ordinary double point. 

The long exact cohomology sequence of the pair $(X,U)$ gives $H^i(X,\Q) = H^i(U,\Q)$ for $i=0,1$ and an exact sequence of mixed Hodge structures
\begin{equation} \label{exact}
0 \to H^2(X,\Q) \to H^2(U,\Q) \to H^3_{\calD}(X,\Q) \to H^3(X,\Q) \to H^3(U,\Q) 
\end{equation}
from which we conclude that $H^k(X,\Q)$ and $H^k(U,\Q)$ are pure of weight $k$ for $k=0,1,2$. Because at all points of $U$ we have at most quotient singularities, $U$ is a rational homology manifold, hence Poincar\'e duality holds for it with rational coefficients: 
$$\mathrm{Hom}(H^k(U,\Q),\Q(-3)) \simeq H^{6-k}_c(U,\Q) \simeq H^{6-k}(X,\Q)$$
This implies the purity of $H^k(X,\Q)$ for $k=4,5$, and $\mathrm{def}(X) = b_2(U)-b_2(X)$. Moreover $W_2H^3(X)$ is the kernel of the restriction map to $H^3(U,\Q)$,  for $H^4_{\calD}(X,\Q)$ is pure of type $(2,2)$ by Poincar\'e duality. Therefore the theorem follows from sequence (\ref{exact}). 
\section{Extensions of mixed Hodge structures and Abel-Jacobi map}
Suppose that $X$ is an ADT-threefold. Let $(\hat{X},E) \to (X,\Sigma)$ be the blowing up. Then $E \to \Sigma$ is a conic bundle with degenerate fibers over the pinch points and triple points of $\Sigma$. One has the exact sequence
\begin{equation}
0 \to H_3(E) \to H_3(\hat{X}) \to H_3(X) \to H_2(E) \stackrel{v}{\to} H_2(\hat{X})\oplus H_2(\Sigma) \label{MV}
\end{equation}
Hence we have the extension of mixed Hodge structures
$$
0 \to W_{-3}H_3(X) \to H_3(X) \to \ker(v) \to 0
$$
and $\ker(v)$ is purely of type $(-1,-1)$. The extension is classified by a homomorphism $b: \ker(v) \to J$ where $J$ is the Jacobian of the Hodge structure $W_{-3}H_3(X)$, see \cite[Theorem 3.31]{PS}. 

On the other hand, let $\mathrm{Hom}_1(\hat{X})$ denote the group of algebraic 1-cycles on $\hat{X}$ which are homologous to zero. We have the Abel-Jacobi map 
$$a: \mathrm{Hom}_1(\hat{X}) \to J'$$
where $J'$ is the intermediate Jacobian of $\hat{X}$. It is equipped with a surjection $p:J'\to J$ with kernel the Jacobian $J(\tilde{\Sigma})$ of the normalisation of $\Sigma$.  
We also have a natural homomorphism $i_\ast: \ker(v) \to \mathrm{Hom}_1(\hat{X})$ where $i:E \to \hat{X}$ is the inclusion.  These fit  into the following diagram 
$$
\begin{array}{ccccccccc}
0 & \to & J(\tilde{\Sigma}) & \to & J'& \stackrel{p}{\to} & J & \to & 0 \\
&&&&\uparrow a & & \uparrow b && \\
&&&& \mathrm{Hom}_1(\hat{X}) & \stackrel{i_\ast}{\leftarrow} & \ker(v) &&
\end{array}
$$ 
\begin{theorem} \label{extensionclass}
$$b = p \circ a \circ i_\ast$$
\end{theorem}
\emph{Proof } The mixed Hodge structure on the (co)homology of $X$ is obtained from the cubical hyperresolution 
$$ \begin{array}{rcl} E & \stackrel{i}{\to} & \hat{X} \\
c\downarrow & & \downarrow\sigma \\
\Sigma & \hookrightarrow & X 
\end{array} $$
and $X$ is homotopy equivalent to the geometric realization 
$$ T=|X_\cdot| := \hat{X} \cup \Sigma \cup E\times[0,1] / \sim$$
where for $y\in E$ we identify  $(y,0)$ with $i(y)\in \hat{X}$ and $(y,1)$ with $c(y) \in \Sigma$. 

Let $T_1$ be the image of $\hat{X}$ in $T$ and let $T_2$ be the image in $T$ of $E \times[0,1]$. Then $T_1$ is isomorphic to $\hat{X}$ and $T_2$ is the mapping cone of $c$. Hence $T_2$ is homotopy equivalent to $\Sigma$. Moreover $T_1\cup T_2 = T$ and $T_1\cap T_2 = E\times (0,1)$ is homotopy equivalent to $E$. Also $(T;T_1,T_2)$ is an excisive triad (see \cite[Sect. III.8]{Dold}), and from its Mayer-Vietoris exact sequence we obtain the sequence (\ref{MV}).

The Jacobian $J'$ can be constructed as follows. Consider $$V:= F^2H^3(\hat{X},\C).$$
Then we have a natural map $h: H_3(\hat{X})  \to \mathrm{Hom}(V,\C)$ given by 
$$h(\delta)(\omega) = \int_{\delta} \omega$$
whose image is a lattice, and 
$$J'= \mathrm{Hom}(V,\C)/h(H_3(\hat{X})).$$ 

The map $a$ is constructed as follows. For $Z \in \mathrm{Hom}_1(\hat{X}) $ choose a singular 3-chain $\Gamma$ on $\hat{X}$ with boundary $Z$. Then 
$$ a(Z): \omega \mapsto \int_{\Gamma}\omega$$
is well-defined up to integrals over elements of $H_3(\hat{X})$. 

The classifying homomorphism $\beta: \ker(v) \to J$ is constructed as follows. Choose a section $s_{\Z}$ of $H_3(X,\C) \to \ker(v_{\C})$ with image in $H_3(X,\Z)$ and another such section $s_F$ with image in $F^{-1}H_3(X,\C)$ (note that $\ker(v)$ is purely of type $(-1,-1)$). Then $\beta$ is the map $s_{\Z}-s_F$ followed by the canonical map $W_{-3}H_3(X,\C) \to W_{-3}H_3(X,\C)/[W_{-3}H_3(X,\Q)\cap H_3(X,\Z) + F^{-1}H_3(X,\C)]$. Note that we may omit $s_F$. 

A section $s_{\Z}$ may be constructed as follows. Represent $\gamma \in \ker(v)$ by an algebraic 2-cycle on $T_1\cap T_2$ and choose 3-chains $\Gamma$ and $\Gamma'$ on $T_1$ and $T_2$ respectively such $\partial \Gamma = \partial \Gamma' = \gamma$. Then $\Gamma-\Gamma'$ is a 3-cycle on $T$ representing $s_{\Z}(\gamma)$. It follows that $$b = p \circ a \circ i_\ast$$
if and only if we can show that $\int_{\Gamma'}v = 0$ for all $v \in F^2H^3(X,\C)$. 

An element of $F^2H^3(X,\C)$ is represented by a closed form $\omega$ of type $(2,1)$ on $\hat{X}$ whose restriction to $E$ is exact. Note that $H^3(E) \simeq H^3(\tilde{E})$ where $\tilde{E}$ is the normalization. Let $j:E \hookrightarrow \hat{X}$ be the inclusion. We can apply the $\partial\bar{\partial}$-lemma   and write $j^\ast\omega = d\eta$ with $\eta$ of type $(2,0)$ on $E$. Then 
$$\int_{\Gamma'}\omega = \int_{\Gamma'}d\eta = \int_{\partial\Gamma'}\eta = \int_\gamma \eta = 0$$
because the restriction of any form of type $(2,0)$ to an algebraic one-cycle is identically zero. 

\section{Ordinary double solids}
We consider ordinary double solids $X \to \PP^3$, with ramification surface $S$ and singular locus $\Sigma$. The \emph{degree of $X$} is defined as the degree of $S$; it is always even. 
In addition to Theorem \ref{ADT} we have the following result:
\begin{theorem} \label{euler-betti}
Let $X$ be an ordinary double solid ramified over a surface $S$ of degree $2d$. Then 
\begin{enumerate}
\item $b_1(X)=b_5(X)=0,\ b_2(X)=1,\ b_4(X)=1+\mathrm{def}(X)$;
\item the Euler number $e(X)$ of $X$ is equal to $8-e(S)$;
\end{enumerate}
\end{theorem}
\paragraph{Proof} The double solid $X$ can be considered as a hypersurface in the weighted projective fourfold $\PP = \PP(1,1,1,1,d)$ given by the equation $X_4^2=F(X_0,X_1,X_2,X_3)$. Note that the only singular point of this weighted projective space is a cyclic quotient singularity. Therefore Poincar\'e duality holds for $\PP \setminus X$. Let $i:X \hookrightarrow \PP$ denote the inclusion. By Lefschetz' hyperplane theorem $i^\ast:H^k(\PP,\Q) \to H^k(X,\Q)$ is an isomorphism for $i \leq 2$. This implies the statements about the Betti numbers. The statement about the Euler number is a general fact about double covers.

\vspace{3mm}\noindent
In the rest of this section we focus on ordinary double solids whose ramification surface $S$ has degree four. There are the following 11 possibilities: six cases where $S$ is irreducible and five cases where $S$ is reducible. We let $\Sigma$ denote the singular locus of $S$. 
\begin{enumerate}
\item $S$ smooth, so $\Sigma = \emptyset$
\item $\Sigma$ is a line ($X$ is called a \emph{spider double solid})
\item $\Sigma$ is a smooth plane conic ($X$ is called a \emph{cyclide double solid})
\item $\Sigma$ is a rational normal cubic ($X$ is called a \emph{twisted cylindrical double solid})
\item $\Sigma$ consists of two skew lines ($X$ is called a \emph{viaduct double solid})
\item $\Sigma$ consists if three concurrent lines which are not in one plane ($X$ is called the \emph{three-axed double solid})
\item $S$ is the union of a smooth cubic surface and a plane intersecting it transversely
\item $S$ is the union of two smooth quadrics in general position
\item $S$ is the union of a smooth quadric and two planes in general position
\item $S$ is the union of a cubic surface with a double line and a plane in general position
\item $S$ is the union of four planes in general position 
\end{enumerate}
The smooth case has been investigated in great detail by Clemens \cite{DS}, who also admits ordinary double points. It follows from Theorems \ref{ADT} and \ref{euler-betti} that all mixed Hodge numers of $X$ are determined as soon as one knows the Euler number, the defect and the number of singular points of type $D$, which is equal to the number of pinch points on $S$. The results are given in the tabel below. The number of singular points of type $D$ is given by the formula 
$$
\sharp(D) -2\sharp(T) = (2\deg S - 8)\deg \Sigma + 4\chi(\calO_\Sigma)
$$
(see \cite{JJ}). 
For the computations of the defects we refer to the forthcoming thesis of the first author. For the defect of the cyclide double solid see Corollary \ref{defect-cyclide}. 

\begin{center}
\begin{tabular}{|c|c|c|c|c|c|}
\hline
Type & Euler number & $\sharp \calD$ & defect & $\dim W_2H^3(X,\Q)$ & $\dim F^2H^3(X,\C)$ \\
\hline
smooth & -16 & 0 & 0 & 0 & 10 \\
\hline
spider & -6 & 4 & 0 & 4 & 3 \\
\hline
cyclide & -2 & 4 & 1 & 3 & 2 \\
\hline
twisted cylindrical & 2 & 4 & 1 & 3 & 0 \\
\hline
viaduct & 4 & 8 & 4 & 4 & 0 \\
\hline
three-axed & 4 & 6 & 3 & 3 & 0 \\
\hline
7 & -4 & 0 & 0 & 0 & 4 \\
\hline
8 & 0 & 0 & 0 & 0 & 2 \\
\hline
9 & 2 & 0 & 0 & 0 & 1 \\
\hline 
10 & 2 & 2 & 0 & 2 & 0 \\
\hline
11 & 4 & 0 & 0 & 0 & 0 \\
\hline
\end{tabular}
\end{center}
Observe that $H^3(X)$ has a pure Hodge structure except in the spider and cyclide cases. 
\section{Geometry of the cyclide double solid} (See \cite[Sect.~11.3]{CAG}.) 
Let $S \subset \PP^3$ be a quartic surface whose singular locus is a plane conic $\Sigma$. The double cover of $\PP^3$ ramified over $S$ is called a \emph{cyclide double solid}.
Its singularities have type $A$ except for the four pinch points of $S$ where $X$ has singularities of type $D$. We let $V$ denote the plane in $\PP^3$ spanned by $\Sigma$. 

Let $f:\bar{Q} \to \PP^3$ be the blowing up of $\PP^3$ with center $\Sigma$. The strict transform $\bar{V}$ of $V$ in $\bar{Q}$ is a plane which can be blown down to a point; the resulting variety $Q$ is isomorphic to a smooth quadric in $\PP^4$ and the birational map $Q \to \PP^3$ is the projection from the point $x \in Q$ corresponding to the blown down plane $\bar{V}$. The exceptional divisor $\bar{E}$ of $f$ is a $\PP^1$-bundle over $\Sigma$ isomorphic to a Hirzebruch surface $F_2$ , and $\bar{V} \cap \bar{E}$ is the section with self-intersection number $-2$ on it. The image $E$ of $\bar{E}$ under $\bar{\pi}_x$ is the intersection of $Q$ with the embedded tangent hyperplane to $Q$ at $x$, and $E$ is a quadric cone. 
The surface $S$ is transformed to its normalization $\tilde{S}$ in $\bar{Q}$ and in $Q$. 

The modifications $\bar{Q}$ and $Q$ of $\PP^3$ induce modifications $\bar{Y}$ and $Y$ of its double cover $X$. We have the following diagram 
$$
\begin{array}{ccccccc}
p_1,p_2,\tilde{S},D & \hookrightarrow & Y & \stackrel{2:1}{\longrightarrow } & Q & \hookleftarrow & \tilde{S},p,E \\
\uparrow &&\uparrow &&\uparrow \bar{\pi}_x && \uparrow \\
V_1,V_2,\tilde{S},\bar{D} &\hookrightarrow & \bar{Y} & \stackrel{2:1}{\longrightarrow } &\bar{Q} & \hookleftarrow &\tilde{S},V,\bar{E}\\
\downarrow  && \downarrow  && \downarrow  && \downarrow \\
V_1,V_2,S,\Sigma & \hookrightarrow & X & \stackrel{2:1}{\longrightarrow } &\PP^3 & \hookleftarrow & S,V,\Sigma
\end{array}
$$
Here $\tilde{S}$ is the smooth intersection of $Q$ with another quadric hypersurface in $\PP^4$ (a Del Pezzo surface of degree four), and $Y \to Q$ has ramification locus $\tilde{S}$. Therefore $Y$ is the intersection of two quadrics in $\PP^5$. The planes $V_1,V_2$ form the preimage of $V$ and are blown down in $Y$ to the two points above $p\in Q$. The surface $\bar{D}$ is a conic bundle over $\Sigma$, whose singular members lie 
over the pinch points on $\Sigma$. It has two sections $\sigma_i=V_i\cap \bar{D},\ i=1,2$ mapping to $V\cap E$, both with self-intersection $-2$, and these are blown down to the points $p_1,p_2$ on $D$, which are ordinary double points. Finally, $\tilde{S} \cap \bar{D} \simeq \tilde{S}\cap D$ is a bisection of the conic bundle, ramified over the pinch points of $\Sigma$. 
\paragraph{Moduli of the cyclide  double solid}
The threefold $Y$ is the smooth intersection of two quadrics in $\PP^5$. The pencil $\Lambda$ of quadrics in $\PP^5$ containing $Y$ has six degenerate members, each of which has a single singularity. All quadrics of the pencil can be simultaneously diagonalized. Therefore we have an action of the group of diagonal matrices
with all entries on the diagonal equal to $\pm 1$ on $Y$.  It follows that $\mathrm{Aut}(Y)$ contains a subgroup $C_2^5$ of order 32, and for  generic $Y$ these are all automorphisms. The subgroup $C_2^5$ contains six involutions whose fixed point set $\tilde{S}$ on $Y$ is a divisor. The double cover $Y \to Q$ is obtained as the quotient of $Y$ by one of these involutions. The isomorphism type of $Y$ is determined by the moduli of the  unordered sextuple of points on $\Lambda$ corresponding to the six degenerate members. The involution chosen corresponds to the choice of one of these points on $\Lambda$. 

The two points $p_1,p_2$ on $Y$ are conjugate by this involution. Hence $\bar{Y}$ is determined by the choice of a pair of points on $Y$, conjugate by the involution chosen. Conversely, the exceptional planes $V_1$ and $V_2$ on $\bar{Y}$ are the only two exceptional planes, so $\bar{Y} \to Y$ is canonically determined by $\bar{Y}$. We conclude
\begin{theorem}\label{moduli}
$Y$ is determined up to isomorphism by a sextuple of points in $\PP^1$; $\bar{Y}$ is obtained from $Y$ by the choice of one these six points plus an unordered pair of  points $\{p_1,p_2\}$ on $Y$ conjugate under the corresponding involution. The surface $\tilde{S}$ is the fix point set of the involution on $Y$ and $S$ is the image of $\tilde{S}$ under the projection from the line connecting the pair of points.  

The sextuple of points determines a pencil of quadrics in $\PP^5$ whose singular members correspond to the six given points, and the choice of one of these quadrics determines the involution, which fixes the singular point of this quadric and the hyperplane spanned by the remaining five singular points. The choice of three smooth quadrics in the pencli then determines the pair $\{p_1,p_2\}$. 
\end{theorem}

\paragraph{The intermediate Jacobian} The intermediate Jacobian of $X$ is the Jacobian of the mixed Hodge structure $H_3(X)$. This mixed Hodge structure has as graded quotients the part of weight $-3$, which is isomorphic to $H_3(\bar{Y}) \simeq H_3(Y)$ and its quotient of weight $-2$, isomorphic to $\ker(v)$ where $v: H_2(\bar{D}) \to H_2(\bar{Y})
\oplus H_2(\Sigma)$ is the natural map (see Sect.~3). We will consider this mixed Hodge structure with its graded polarization, which is the intersection form on $H_3(Y)$ for the weight $-3$
part and the intersection form on $H_2(\bar{D})$ restricted to $\ker(v)$ on the weight $-2$
quotient. 

\begin{lemma} $\ker(v)$ with its intersection form is a lattice of type $A_3$.
\end{lemma}
{\em Proof } The homology of the conic bundle $\bar{D} \to \Sigma$ is free abelian on generators the homology classes of (i) a section $\sigma_1$, (ii) a smooth fibre, (iii) those components of reducible fibres which do not intersect the section $\sigma_1$. Indeed, this is clear by considering the case of a $\PP^1$-bundle and the behaviour of homology under blowing-up. Let $f$ denote the class of a fibre, and let $a_1,a_2,a_3,a_4$ denote the classes of the four components of the reducible fibres over the pinch points on $\Sigma$ which do not intersect $\sigma_1$. The group $H_2(\bar{Y})$ is free of rank 3, a  basis of its dual $H^2(\bar{Y})$ is $\{h,v_1,v_2\}$ with $h$ the pull-back of a positive generator of $H^2(X)$ and $v_i=[V_i]$. The pull-back of these generators to $H^2(\bar{D})$ are $2f,\sigma_1,\sigma_2$, and the image of $H^2(\Sigma)\to H^2(\bar{D})$ is spanned by $f$. Hence $\ker(v)$ is the orthogonal complement in the lattice $H_2(\bar{D})$ of the sublattice spanned by $\sigma_1,\sigma_2,f$. This is generated by the differences $a_i-a_j$. Taking into account that $a_i^2=-1$ this gives a lattice of type $A_3$.
\begin{corollary} \label{defect-cyclide}
The cyclide double solid has defect 1.
\end{corollary}

\section{Torelli theorem for the cyclide double solid}
 We finally show which of the data in Theorem \ref{moduli} can be recovered from the polarized intermediate Jacobian of $X$. In the first place, the abelian part $J$ of this intermediate jacobian is isomorphic to the polarized intermediate Jacobian of $Y$, which in turn is isomorphic to the polarized Jacobian of the hyperelliptic curve $C$ which is ramified over $\Lambda$ over the six exceptional points of the pencil. The Torelli theorem for curves tells us that we can recover $C$ from its polarized jacobian, and hence we recover the threefold $Y$ from it. 

Much is known about the geometry of $Y$; by identifying one of the smooth quadrics containing $Y$ with the Grassmannian of lines in $\PP^3$ we get an isomorphism of $Y$ with a \emph{quadric line complex}, see \cite[Ch.~6]{GH}. The variety $F$ of lines on $Y$ is a homogeneous space under the abelian surface $J$. The isomorphism type of the intermediate Jacobian of $X$ is determined by a classifying homomorphism $b:\ker(v) \to J$ strongly related to the Abel-Jacobi map, as described in Section 3. Recall that $\ker(v)$ is a lattice of type $A_3$. 

Recall the group $ \mathrm{Hom}_1(Y)$ of rational equivalence classes of algebraic 1-cycles on $Y$ which are homologous to zero. The Abel-Jacobi map 
$$a:\mathrm{Hom}_1(Y) \to J$$ 
is an isomorphism. For a fixed line $\ell_0 \in F$ on $Y$, the map $\ell \to a(\ell-\ell_0)$ defines an isomorphism between $F$ and $J$. 

We first note that the map $b:\ker(v) \to J$ is injective. Choose a basis of $\ker(v)$ which is a root basis for the root system of type $A_3$. The images of these roots are three points on $J$. 

The surface $F$ is covered by curves $B_\ell \simeq C$ for $\ell \in F$ as follows: 
$$B_\ell = \mbox{ closure of }\{m\in F\setminus \{\ell\}\mid m\cap \ell \neq \emptyset\}\ .$$
Because $B_\ell$ has genus 2 and $F$ has trivial canonical bundle, the adjunction formula shows that $B_\ell\cdot B_\ell = 2$. Hence $B_\ell \neq B_{\ell'} \Rightarrow \sharp(B_\ell \cap B_{\ell'}) \leq 2$. 

If $\ell_0$ is one of the sixteen lines on the Del Pezzo surface $\tilde{S}$, then $B_\ell = B_{\ell_0} - (\ell - \ell_0)$. 
For $\ell \in F$ we may also define a curve $C_\ell \subset J$ by 
$$C_\ell = \{(m)-(\ell)\mid m\in B_\ell\} = B_\ell - \ell$$
which is isomorphic to $B_\ell$. See \cite[p.~783]{GH}. 
\paragraph{Claim:} 
$$C_\ell = C_{\ell'} \Longleftrightarrow (\ell)-(\ell') \in J[2]$$
where $J[2]$ is the kernel of multiplication by 2 on $J$. 
\paragraph{Proof of the claim. } Suppose that $C_\ell=C_{\ell'}$. Then 
$$B_{\ell_0} - (\ell - \ell_0) -\ell = B_{\ell_0}-(\ell'- \ell_0) - \ell'$$
so $2(\ell-\ell') = 0$ in $A$. 

\vspace{3mm}\noindent
Note that the images of the three roots all belong to one and the same $C_\ell$ (they are of the form $a_2-a_1,a_3-a_1,a_4-a_1$ and choose $\ell=a_1$). Hence $a_1\in F$ is determined by them up to translation by an element of order two in $J$. Fix a choice for $a_1$. Then $a_2,a_3,a_4$ are also determined. They intersect on $Y$ in a point $p_1$. The choice of an involution $g$ of $Y$ with a surface of fixed points then determines the point $p_2=g(p_1)$. 

A different choice for $a_1$, say $a_1^\prime$, is related to $a_1$ by translation by a point of order two on $J$. The action of $J[2]\simeq C_2^4$ on $F$ however is induced by the action of $\mathrm{SAut}(Y) \simeq C_2^4$ on it where $\mathrm{SAut}(Y)= \mathrm{Aut}(Y) \cap \mathrm{PSL}(6,\C)$. Indeed, an automorphism of $Y$ with determinant 1 has on $\PP^5$ a fixed line (not on $Y$) and a fixed $\PP^3$ which intersects $Y$ in a smooth curve of genus one. Hence such an automorphism induces an involution without fixed points on $F$ and therefore must be a translation by a point of order two of $J$. 

We conclude:
\begin{theorem}
Given the polarized mixed Hodge structure on $H_3(X)$, where $X$ is a cyclide ordinary double solid, there are six possibilities for $X$ itself. These possibilities correspond to the choice of a singular quadric containing $Y$. 
\end{theorem}

\paragraph{Remark.} There is a 1-1 correspondence between isomorphism classes of cyclide double solids (or, equivalently, of cyclide surfaces)
and 
moduli of $5+1+3$ points on a projective line. The five points determine the Del Pezzo surface, the extra point determines the smooth quadric $Q$ contianing it and its double cover, and the 1+3 points determine the branch points of the double cover of $\Sigma$. 

This is seen most clearly in $\PP^4$. The situation is uniquely determined there by a pencil of quadrics and a point $q$ outside the intersection of these quadrics. The point $q$
lies on a unique quadric $Q$ of
the pencil, corresponding to a sixth member of the pencil. The tangent hyperplane $T$ to $Q$ at $q$ is tangent to three other quadrics of the pencil;
these correspond to the  other three points. 

As we just have shown, the mixed Hodge structure on $H_3(X)$ also determines 9 points on a line, but arranged in a $6+3$ configuration. The passage from $X$ to the polarized mixed Hodge structure on $H_3(X)$ is the functor which forgets the choice of one of these six points.

\end{document}